\documentclass{svmult}

\usepackage{makeidx}
\usepackage{multicol}

\usepackage{amsmath}
\usepackage{amssymb}

\makeindex

\def\semi {\times_s}
\def\Chi{{\chi}}
\def\Z{{\mathbb{Z}}}
\def\R{{\mathbb{R}}}
\def\N{{\mathbb{N}}}
\def\T{{\mathbb{T}}}
\def\={\ = \ }

\def\supp{{\mathrm{supp}}}

\def\w_hat{\hat w}
\def\f_hat{\hat f}
\def\psih{\hat \psi}

\def\fhat{\f_hat}
\def\psihat{\psih}

\def\L^2(\rnhat){L^2(\widehat{\R}^n)}
\def\rnhat{\widehat{\R^n}}

\def\dim#1{{\rm dim}_#1}
\def\myvec#1{\vec{#1}}

\def\vecgamma{\myvec{\gamma}}

\def\vecv{\myvec{v}}
\def\vecvone{\myvec{v}_1}

\def\vecvtwo{\myvec{v}_2}
\def\vecvthree{\myvec{v}_3}
\def\vecvfour{\myvec{v}_4}
\def\vecvfive{\myvec{v}_5}
\def\vecvn{\myvec{v}_n}
\def\calS {\mathcal{S}}
\def\calT {\mathcal{T}}

\def\bysame{----------------}

\def\dt{\; dt}

\begin{document}

\title{Explicit cross sections of singly generated group actions}
\author{David Larson \thanks{Partially supported by NSF grant DMS-0139386} \inst{1} \and Eckart Schulz \thanks{Supported by a research grant from
Suranaree University of Technology}
\inst{2}
\and Darrin Speegle \inst{3} \and Keith F.~Taylor \inst{4}}

\institute{Department of Mathematics,
Texas A\&M University,
College Station, TX 77843-3368
\texttt{larson@math.tamu.edu} \and 
School of Mathematics, Suranaree University of Technology, 
111 University Avenue, Nakhon Ratchasima, 30000, Thailand
\texttt{eckart@math.sut.ac.th} \and
Department of Mathematics, 
Saint Louis University, 
221 N Grand Blvd, 
St. Louis, MO 63103
\texttt{speegled@slu.edu} \and
Department of Mathematics and Statistics,
Dalhousie University,
Halifax, Nova Scotia,
Canada B3H 3J5
\texttt{keith.f.taylor@dal.ca}}

\maketitle
\begin{abstract}
We consider two classes of actions on $\mathbb{R}^n$ - one continuous and one discrete.  For matrices of the form $A = e^B$ with $B \in M_n(\R)$, we
consider the action given by $\gamma \to \gamma A^t$.  We characterize the matrices $A$ for which there is a cross-section for this action.  The discrete
action we consider is given by
$\gamma
\to
\gamma A^k$, where $A\in GL_n(\R)$.  We characterize the matrices $A$ for which there exists a cross-section for this action as well.  We also
characterize those $A$ for which there exist special types of cross-sections; namely, bounded cross-sections and finite measure cross-sections.  Explicit
examples of cross-sections are provided for each of the cases in which cross-sections exist.  Finally, these explicit cross-sections are used to
characterize those matrices for which there exist MSF wavelets with infinitely many wavelet functions.  Along the way, we generalize a well-known
aspect of the theory of shift-invariant spaces to shift-invariant spaces with infinitely many generators.
\end{abstract}
\section{Introduction}

In discrete wavelet\index{wavelet!discrete}\index{discrete wavelet} analysis on the line, the classical approach is to dilate and translate a single
function, or \emph {wavelet}, so that the resulting system is an orthonormal basis for $L^2(\R)$.  More precisely,  a wavelet is a function $\psi \in
L^2(\R)$ such that 
\[\{2^{j/2} \psi(2^jx + k) : k,j\in \Z\}
\]
forms an orthonormal basis of $L^2(\R)$.  

In multidimensional discrete wavelet analysis, the approach is similar.  Fix a matrix $A\in GL_n(\R)$ and a full rank lattice $\Gamma$.
A collection of functions $\{\psi^i: i=1,\ldots,N\}$ is called an \emph {$(A, \Gamma)$ orthonormal wavelet of order $N$} if 
dilations by $A$ and translations by $\Gamma$,
\[ \{|\det A|^{j/2} \psi^i(A^jx + k): i=1,\dots,N, j\in \Z, k\in \Gamma\},
\] 
forms an orthonormal 
basis for $L^2(\R^n)$.  In this generality, there is no characterization (in
terms of $A$ and $\Gamma$) of when wavelets exist.  It was shown in \cite{DLS97} that, 
if $A$ is expansive (that is, a matrix whose eigenvalues all have
modulus greater than 1) then there does exist an orthonormal wavelet.  
A complete characterization of such wavelets in terms of the Fourier transform
was given in \cite{HLW02}.  The non-expansive case remains problematic.

It is also possible to study the continuous version of wavelet analysis.  \index{wavelet!continuous} \index{continuous wavelet}
Consider the full affine group of motions given by 
$GL_n(\R) \times \R^n$ with multiplication given by $(a,b) (c,d) = (ac, c^{-1}b + d)$.  
We are interested in subgroups of the full affine group of motions of the form
\[
	G = \{(a,b) : a\in D, b\in \R^n\},
\]
where $D$ is a subgroup of $GL_n(\R)$.  In this case, $G$ is the 
semi-direct product $D\semi \R^n$.  Now, if we define the unitary operator 
$T_g$ for $g\in G$ by
\[
	\bigl(T_g\psi\bigr)(x)  = |\det a|^{-1/2}\psi(g^{-1}(x)), 
\]
then the continuous wavelet transform is given by
\[
	\langle f, \psi_g \rangle := \int_{\R^n} f(x) \overline{(T_g\psi)(x)} \, dx, 
\]
which is, of course, a function on $G$. The function $\psi$ is a $D$-continuous 
wavelet if it is possible to reconstruct all functions $f$ in $L^2(\R^n)$
via the following reconstruction formula:
\[
	f(x) = \int_G \langle f, \psi_g \rangle \psi_g(x)\, d\lambda(g),
\]
where $\lambda$ is Haar measure on $G$.

There is a simple characterization of continuous wavelets, given in \cite{WW01}.  
\begin{theorem}\label{ctswaveletchar}\cite{WW01}  Let $G$ be a subgroup of 
the full affine group of the form $D\semi \R^n$.  A function $\psi \in L^2(\R^n)$ is a
$D$-continuous wavelet if and only if the Calder\'on condition \index{Calder\'on condition}\index{continuous wavelet!characterization}
\begin{equation}\label{calderon1}
	\int_D |\hat\psi(\xi a)|^2 \, d\mu(a) = 1 \hskip .1in a.e. \ \xi \text{ in } \rnhat
\end{equation}
holds, where $\mu$ is left Haar measure for $D$.
\end{theorem}

In this paper, we will always assume one of the two following
cases, which for the purposes of this paper will be the {\it singly generated subgroups} of $GL_n(\R)$.  \index{group!singly generated}
\begin{enumerate}
\item\label{one} $D = \{A^k: k\in \Z\}$ for some $A\in GL_n(\R)$, or
\item\label{two} $D = \{A^t: t\in \R\}$ for some $A = e^B$, where $B\in M_n(\R)$.
\end{enumerate}
We will say that $D$ is generated by the matrix $A$.
Applying Theorem \ref{ctswaveletchar} to these cases gives the following characterizations. 
\begin{proposition}\label{prop1}
\begin{enumerate}
\item Let $A \in GL_n(\R)$ and denote the dilation group $D = \{A^k: k\in \Z\}$.  Then, $\psi \in L^2(\R^n)$  is a $D$-continuous
wavelet if and only if 
                \[       \sum_{k \in \Z} | \psihat(\xi A^k)|^2 = 
1         \]
        for almost all $\xi \in \rnhat$.
\item Let $D = \{A^t: t\in \R\}$ for some $A = e^B$, where $B\in M_n(\R)$.  Then $\psi \in L^2(\R^n)$  is a 
$D$-continuous wavelet if and only if
                \[      \int_{\R} | \psihat(\xi A^t)|^2 \dt = 
1    \]
        for almost all $\xi \in \rnhat$.
\end{enumerate}
\end{proposition}

In the case that $D$ is generated by a single matrix as above, a complete characterization of 
matrices for which there exists a continuous wavelet is
given in \cite{LWWW02}.
\begin{theorem}\label{charaddmisiblegroups} Consider the dilation group $D$ as in case \ref{one} or \ref{two} above.  There exists a continuous wavelet
if and only if $|\det(A)| \not = 1$.
\end{theorem}
The wavelets constructed in \cite{LWWW02} are of the form $\hat\psi = \Chi_K$, for some set $K$.   
One drawback to the proof in \cite{LWWW02}
is that, while the proof is constructive, the sets $K$ that are constructed are written as the countable union of set differences of sets consisting
of those points whose orbits land in a prescribed closed ball a positive, finite number of times.  Hence, it is not clear whether the set constructed in
the end can be chosen to be ``nice" or  easily described.  

The purpose of this article is two-fold.  First, we will give explicitly defined, easily verified sets $K$ such that $\Chi_K$ is the Fourier transform
of a continuous wavelet.  Here, we will exploit the fact that we are in the singly generated group case to a very large extent.  We will
also obtain a characterization of matrices such that the set $K$ can be chosen to be bounded as well as a characterization of matrices such that the
only sets $K$ that satisfy (\ref{calderon1}) have infinite measure.

Second, we will show how to use these explicit forms to characterize those matrices 
such that there exists a discrete wavelet of order infinity. 
Note that this seems to be a true application of the form of the sets $K$ in 
section \ref{Cross-sections}, as it is not clear to the authors how to use
the proof in \cite{LWWW02} (or the related proof in \cite{ILP98}) to achieve the same result.

\section{Cross-sections}\label{Cross-sections}

Throughout this section, we will use vector notation to denote elements of $\rnhat$,
and $m$ will denote the Lebesgue measure on $\R^n$.  Multiplication of a vector with a matrix will be given by $\vecgamma A$, and we will reserve the
notation $A^t$ as ``$A$ raised to the $t$-power".  In the few places we need the transpose of a matrix, we will give it a separate name.

\begin{definition}\label{cross-section}\index{cross-section}
         A Borel set $S \subset \rnhat$ is called a 
        \emph{cross-section} for the continuous action $\vecgamma \to 
\vecgamma A^t \ (t \in \R)$ if
\begin{enumerate}
        \item $\bigcup_{t \in \R} SA^t = \rnhat \backslash N$ 
for some set $N$ of measure zero and
        \item $SA^{t_1} \cap SA^{t_2} = \emptyset$ whenever $t_1 \neq 
t_2 \in \R$.
\end{enumerate}
Similarly,  a Borel set $S \subset \rnhat$ is called a 
        \emph{cross-section} for the discrete action $\vecgamma \to 
\vecgamma A^k \ (k \in \Z)$ if
\begin{enumerate}
        \item $\bigcup_{k \in \Z} SA^k = \rnhat \backslash N$ 
for some set $N$ of measure zero and
        \item $SA^{j} \cap SA^{k} = \emptyset$ whenever $j \neq k \in 
\Z$.
\end{enumerate}
\end{definition}

Note that we have defined cross-sections using left products, which will eliminate the need for taking transposes in section \ref{sect13}.

Note also that if $S$ is a cross-section for the continuous action, then 
$\{ \vecgamma A^t: \vecgamma \in S, 0 \leq t < 1 \}$
is a cross-section for the discrete action.  Cross-sections are sometimes referred to as multiplicative tiling sets.\index{tiling set!multiplicative}

\begin{remark}\label{changebasis}
Let $S$ be a cross-section for the action $\vecgamma \to \vecgamma A^k$.  
Then, $SJ^{-1}$ is a cross-section for the action $\vecgamma \to \vecgamma JA^kJ^{-1}$, and similarly,
for the continuous action $\vecgamma \to \vecgamma A^t$, 
where $A = e^B$, $SJ^{-1}$ is a cross-section for the continuous action $\vecgamma \to \vecgamma \tilde A^t$, where $\tilde
A = e^{JBJ^{-1}}$. 
\end{remark}

To begin with cross-sections for the continuous action, let $A=e^B 
\in GL_n(\R)$ be given, where by the preceding remark we may assume that $B$ is in real Jordan 
normal form. \index{Jordan form!real}
Then, $B$ is a block diagonal matrix,
and a block corresponding to a real eigenvalue $\alpha_i$ is of the 
form
\[
    B_i = \begin{pmatrix}
                                \alpha_i & 1 && (0) \\
                                & \ddots & \ddots & &\\
                                & &\ddots & 1  \\
                                (0) &&  & \alpha_i
                        \end{pmatrix}
\]
while a block corresponding to a complex pair of eigenvalues  
$\alpha_i \pm i \beta_i$ with $\beta_i \neq 0$ is of the form
\[
                        B_i = \begin{pmatrix}
                                D_i & I_2 && (0) \\
                                & \ddots & \ddots & &\\
                                & &\ddots & I_2  \\
                                (0) &&  & D_i
                        \end{pmatrix}
                \qquad \text{with} \qquad 
                        \begin{matrix}
                                D_i &= \begin{pmatrix}
                                        \alpha_i & \beta_i \\
                                        -\beta_i & \alpha_i
                                        \end{pmatrix} \\ \\
                                I_2 &= \begin{pmatrix}
                                        1 & 0 \\
                                        0 & 1
                                        \end{pmatrix} 
                        \end{matrix}
\]
In this basis, $A^t$ is again a block diagonal matrix,
and its blocks are of the form
        \[ A_i  \=      e^{tB_i} \=   \begin{pmatrix}
                         \lambda_i^t E_{_i}(t)  & t \lambda_i^t 
E_{_i}(t)  & \frac{t^2}{2!} \lambda_i^t E_{_i}(t)  & \dots\dots   & 
                        \frac{t^{m-1}}{(m-1)!}  \lambda_i^t 
E_{_i}(t)  \\[0.2cm]
                         &       \lambda_i^t E_{_i}(t) &  t 
\lambda_i^t E_{_i}(t)  & &  \vdots \\[0.2cm]
                        & \ddots & \ddots & \ddots & \vdots  \\[0.2cm]
                        &  & \ddots &  t \lambda_i^t E_{_i}(t)   & 
\frac{t^2}{2!} \lambda_i^t E_{_i}(t)   \\[0.4cm]
                        &  & & \lambda_i^t E_{_i}(t)    &t 
\lambda_i^t E_{_i}(t)    \\[0.4cm]
                        (0) &&  &&  \lambda_i^t E_{_i}(t) 
                \end{pmatrix}
        \]
        with $\lambda_i=e^{\alpha_i}$ and $E_i(t)=1$ or $ 
E_i(t)=E_{\beta_i}(t)  = 
                \begin{pmatrix}
                                \cos{\beta_i t} & \sin{\beta_i t} \\
                                -\sin{ \beta_i t} & \cos{\beta_i t}
                \end{pmatrix}$
        depending on whether this block corresponds to a real 
eigenvalue or a pair of complex eigenvalues of $B$.
        The eigenvalues of $A$ are thus
                $e^{\alpha_i}$ and $e^{\alpha_i} e^{\pm i \beta_i}$, 
respectively.

        For ease of notation, when referring to a specific block 
$A_i$ of $A$ we will drop the index $i$. 
        Furthermore, $\vecvone,\dots,\vecvn$
        will denote a Jordan basis of $\rnhat$ chosen so that 
this block  under discussion is the first block, and $(x_1,\dots,x_n)$
will denote the components of a vector $\vecgamma$ in this basis.

\begin{theorem}\label{thm1}\index{cross-section!characterization}
Let $A=e^B$, where $B \in M_n(\R)$ is in Jordan-normal form.
There exists a cross-section for the continuous action $\vecgamma \to 
\vecgamma A^t$ if and only if $A$ is not orthogonal.
\end{theorem}
\begin{proof}
Assume that $A$ is not orthogonal. Then
at least one of the following four situations, formulated in terms of 
the eigenvalues of $B$, will always apply. 

\medskip

\emph{Case 1:}  $B$ has a real eigenvalue $\alpha \neq 0$. A 
corresponding block of $A^t$,  
which we may assume to be the first block, is of the form  
        \[    \begin{pmatrix}
                        \lambda^t &  t \lambda^t  & \dots & 
\frac{t^{m-1}}{(m-1)!}\lambda^t \\
                         & \ddots& \ddots & &\\
                        &  & \ddots &  t \lambda^t \\
                        (0) &&  &  \lambda^t
                \end{pmatrix}
        \]
        with $\lambda=e^{\alpha} \neq 1$. Set
        \[ S = \{ \pm \vecvone  \} \times 
span(\vecvtwo,\dots,\vecvn).  \]
        Then $S$ is a cross-section and  $\bigcup_{t \in \R} SA^t = 
                \{ (x_1,\dots,x_n) \in \rnhat : x_1 \neq 0 
\}$. 
\medskip

\emph{Case 2:} $B$ has a complex pair of eigenvalues $\alpha \pm i \beta$ 
with $\alpha\neq 0$, $\beta >0$.
        At least one block of $A^t$ is then of the form
         \begin{equation} %\label{matrix:case2}
                \begin{pmatrix}
                   \lambda^t E_{\beta}(t) & t \lambda^t  E_{\beta}(t) & \dots 
			& \frac{t^{m-1}}{(m-1)!}\lambda^t E_{\beta}(t)\\
                         & \ddots& \ddots & &\\
                        & & \ddots &  t \lambda^t  E_{\beta}(t)  \\
                        (0) &&  & \lambda^t E_{\beta}(t)
                \end{pmatrix},
          \end{equation}
          and replacing $B$ with $-B$ if necessary, we may assume that
	 $\lambda=e^{\alpha} > 1$. One easily checks that
 \[      S = \{ s \vecvone : 1 \leq s < 
  	\lambda^{2\pi/\beta} \} \times span(\vecvthree,\dots,\vecvn )     \]
        is a cross-section and $\bigcup_{t \in \R} \; SA^t=
	\{  (x_1,\dots,x_n) \in \rnhat : x_1^2+x_2^2 \neq 0 \}$.

\medskip

\emph{Case 3:}  $B$ has an eigenvalue $\alpha=0$ and at least one of 
the blocks of $B$ belonging to this eigenvalue has nontrivial 
nilpotent part. Then the corresponding block of $A^t$  is of the form
        \begin{equation} \begin{pmatrix}
                        1 & t && (*) \\
                         & \ddots& \ddots & &\\
                        &  & \ddots & t \\
                        (0) &&  &  1
                \end{pmatrix}           \label{matrix1}
        \end{equation}
        and is of at least size $2 \times 2$.
        We set
                \[  S = \left \{ \; s \vecvone : s \in \R\backslash 
\{0\} \; \right\}  
                        \times span(\vecvthree,\dots,\vecvn     )  \]
        so that $S$ is a cross-section and  $\bigcup_{t \in \R} \; 
SA^t=\{  (x_1,\dots,x_n) \in \rnhat : x_1 \neq 0 \} $.

\medskip

\emph{Case 4:} $B$ has a purely imaginary pair of eigenvalues $\pm 
i\beta$, $\beta>0$, and at least one of 
the blocks of $B$ belonging to this pair has nontrivial nilpotent 
part. Then the corresponding block of $A^t$ is of the form
        \begin{equation}     
                \begin{pmatrix}
                        E_{\beta}(t) &  tE_{\beta}(t) && (*) \\
                         & \ddots& \ddots & &\\
                        & & \ddots & tE_{\beta}(t) \\
                        (0) &&  & E_{\beta}(t)
                \end{pmatrix}           \label{matrix2}
        \end{equation}
        and is of at least size $4 \times 4$. Set
                \[  S =  \{ p\vecvone + q \vecvthree + s \vecvfour: 
	p>0, \  0 \leq q < \frac{2\pi}{\beta} \, p, \ s \in \R \} 
                    \times  span(\vecvfive,\dots,\vecvn ).\]
        Since this is the least intuitive case, let us verify in 
detail that $S$ is a cross-section. For convenience, we group the 
first four coordinates of a vector $\vecgamma \in \rnhat$ 
into two pairs, and write
        \[  \vecgamma = \bigl ( \; (x_1,x_2), \; (x_3,x_4) , \; x_5, 
x_6 \, \dots, x_n \bigr ) , \]
        so that
        \[  \vecgamma A^t \=  \bigl ( \; (x_1,x_2)E_{\beta}(t), \; 
t(x_1,x_2)E_{\beta}(t)+(x_3,x_4)E_{\beta}(t), \dots \bigr )  . \]
        Now $E_{\beta}(t)$ acts by rotation through the angle $\beta 
t$, so whenever $x_1^2+x_2^2\neq 0$
        then there exists $t_1 \in \R$ such that
                \[   (x_1,x_2)E_{\beta}(t_1) \= (p,0) \]
        for some $p >0$. Then
        \[  \vecgamma A^{t_1} \=   (  p,0,  t_1p + y_3, y_4, \dots  
)   \]
        where $(y_3,y_4)=(x_3,x_4)E(t_1)$. So if we set $t_2=t_1+ k 
\, \frac{2 \pi}{\beta}$ for some integer $k$, then
        \[  \vecgamma A^{t_2} \=   (  p,0,  k \frac{2 \pi p}{\beta} 
+  t_1 p + y_3, y_4, \dots  )    .   \]
        Now there exists a $k$ such that
        \[ 0 \leq  k \frac{2 \pi p}{\beta} +  t_1 p + y_3 < \frac{2 
\pi p}{\beta} \]
        and for this choice of $k$, $\vecgamma A^{t_2} \in S$.
 We conclude that
 \[  \bigcup_{t \in \R} \; SA^t \= \{ (x_1,\dots,x_n) \in 
	\rnhat: x_1^2 + x_2^2 \neq 0 \} . \]
        Suppose now that 
                \[ \vecgamma_1 A^{t_1} \= \vecgamma_2 A^{t_2} \]
        for some $\vecgamma_1, \vecgamma_2 \in S, \ t_1, t_2 \in \R$. 
Equivalently,
                \[    \vecgamma_1 \=  \vecgamma_2 A^t   \]
        for some $t$. If $\vecgamma_1=(p_1,0,q_1,s_1,\dots)$ and  
$\vecgamma_2=(p_2,0,q_2,s_2,\dots)$
        then
        \[ \bigl( \; (p_1,0), \; (q_1,s_1), \dots \bigr ) \= \bigl ( 
\; (p_2,0)E_{\beta}(t), \; t(p_2,0)E_{\beta}(t) \, + \, (q_2,s_2) 
E_{\beta}(t), \dots \bigr ) \]
        so that
                \begin{align*}     
                         (p_1,0) &=  (p_2,0)E_{\beta}(t) \\
                        (q_1,s_1) &=  t(p_2,0)E_{\beta}(t) \, + \, 
(q_2,s_2) E_{\beta}(t).
                \end{align*}
        The first equality gives $p_1=p_2$ and $t=\frac{2\pi}{\beta} 
\, k$ for some integer $k$. Then the second equality reads
                \[      (q_1,s_1) =  ( \frac{2\pi}{\beta} \, k p_1 \, 
+ \, q_2, s_2 )  \]
        which gives $s_1=s_2$ and because  $0 \leq q_2, q_1 < 
\frac{2\pi p}{\beta} $, also that $k=0$ and $q_1=q_2$.
        Thus, $S$ is indeed a cross-section.

\medskip

Now suppose to the contrary that $A$ is orthogonal, but there exists
a cross-section $S$. Then
\[      T = \{ \vecgamma A^t : \vecgamma \in S, \ 0 \leq t < 1 \} \]
is a cross-section for the discrete action of $A$ on $\widehat{\R^n}$.
Note that $A$ maps the closed unit ball $B_1(0)$ onto itself,
so if $T_o= T \cap B_1(0)$
then
\[      B_1(0) = \bigcup_{k \in \Z} T_oA^k,     \]
except for a set of measure zero, and this union is disjoint. Then
\[      m\bigl ( B_1(0) \bigr ) = \sum_{k \in \Z} m(T_o A^k)
        = \sum_{k \in \Z} m(T_o) \in \{0, \infty\}
\]
which is impossible.
\end{proof}

\begin{remark}\label{rem3}
The cross-sections constructed in the proof above allow for a change 
of variables 
to integrate along the orbits. 

For example in case 3), given $\vecgamma =(x_1,x_2,\dots,x_n) \in 
\rnhat$ with $x_1 \neq 0$, we set
        \[  \vecgamma = F(t,s,a_3,\dots, a_n) = 
(s,0,a_3,\dots,a_n)A^t   \]
where $s \neq 0$. The Jacobian of this transformation is
\[
        \begin{vmatrix}
                (s,0,a_3,\dots,a_n)B \\
                \vecvone \\
                \vecvthree \\
                        \vdots\\
                \vecvn 
        \end{vmatrix} \; \det(A)^t = \begin{vmatrix}
                0   & s  & * & \dots & * \\
                1   & 0 & 0 & \dots & 0 \\
                0   & 0 & 1 & \dots & 0 \\
                        \vdots \\
                0   & 0 & 0 & \dots & 1 \\
        \end{vmatrix} \; \det(A)^t
                = -s  \delta^t \neq 0,
\]
so that for $\fhat \in \L^2(\rnhat)$,
\[
        \int_{\rnhat} \fhat(\vecgamma) \, d\vecgamma \=  
\int_{\R^{n-2}} \int_{\R\backslash \{0\}} 
\int_{\R}
                 \fhat( \, (s,0,a_3,\dots,a_n)A^t \, ) \, |s| 
\delta^t \, dt \, ds \, da_3 \cdots da_n . 
\]

In case 4), given $\vecgamma =(x_1,x_2,\dots,x_n) \in \rnhat$ 
with $x_1^2 + x_2^2 \neq 0$, we set
        \[  \vecgamma = F(t,p,q,s,a_5,\dots, a_n) = 
(p,0,q,s,a_5,\dots,a_n)A^t  \]
where $p>0$, $ 0 \leq q < 2\pi p /\beta$. The Jacobian of this 
transformation is
\begin{align*}
        \begin{vmatrix}
                (p,0,q,s,a_5,\dots,a_n)B \\
                \vecvone \\
                \vecvthree \\
                        \vdots\\
                \vecvn 
        \end{vmatrix} \; \det(A)^t
        &\=
        \begin{vmatrix}
                0 & \beta p     & * & *         & * & \dots & * \\
                1   & 0         & 0 & 0         & 0 & \dots & 0 \\
                0   & 0         & 1 & 0         & 0 & \dots & 0 \\
                        \vdots \\
                0   & 0         & 0 &           & 0 & \dots & 1 \\
        \end{vmatrix} \; \det(A)^t      
        \\ \\
                &\= - \beta p \delta^t \neq 0
\end{align*}
since $\beta \neq 0$. Thus,
\begin{multline*}
        \int_{\rnhat} \fhat(\vecgamma) \, d\vecgamma \=  
\int_{\R^{n-4}} \int_0^\infty  \int_{\R}
                \int_0^{2\pi p/\beta}  \int_{\R}  
                \fhat( \,  (p,0,q,s,a_5,\dots,a_n)A^t \, ) \, 
|\beta|  p  \delta^t \\
                                \, dt  \, dq \, ds \, dp \, da_5 
\dots da_n .
\end{multline*}
\end{remark}

Any invertible matrix gives rise to a discrete action on 
$\rnhat$, and nearly always there will exist
a cross-section for this action:

\begin{theorem}\label{thm2}\index{cross-section!discrete}
Let $A \in GL_n(\R)$ be in Jordan normal form, and consider the 
discrete action $\vecgamma \to \vecgamma A^k$.
\begin{enumerate}
        \item There exists a cross-section if and only if $A$ is not 
orthogonal.
        \item\label{part2} There exists a cross-section of finite measure if and 
only if $|\det(A)| \neq 1$.\index{cross-section!finite measure}
        \item There exists a bounded cross-section  if and only if 
the (real or complex) eigenvalues of $A$ have all modulus $>1$ or all 
modulus $<1$.\index{cros-section!bounded}
\end{enumerate}
\end{theorem}
\begin{proof}
        To prove the first assertion, choose a Jordan basis 
$\vecvone,\dots,\vecvn$
        so that the Jordan block of $A$ 
        under discussion is the first block. Each Jordan block will 
be an
        upper diagonal matrix of the form 
                \begin{equation*}   
                \begin{pmatrix}
                        \lambda^k E_{\beta}(k) & \binom{k}{1} 
\lambda^{k-1}  E_{\beta}(k-1) & \dots \ \dots & \binom{k}{m-1} 
\lambda^{k-m+1}E_{\beta}(k-m+1) \\
                         & \ddots& \ddots & \vdots  &\\
                        & & \ddots &  \binom{k}{1} \lambda^{k-1}  
E_{\beta}(k-1)  \\[0.2cm]
                        (0) &&  & \lambda^k E_{\beta}(k)  
                \end{pmatrix}
                \end{equation*}
        where $E_{\beta}=1$ if this block corresponds to a real 
eigenvalue $\lambda$,
        and $E_{\beta}$ is a rotation if it belongs to a complex pair 
$\lambda e^{\pm i \beta}$
        of eigenvalues. By a change of basis, we can always simplify 
this block to
        \begin{equation}   
                \begin{pmatrix}
                        \lambda^k E_{\beta}(k) & \binom{k}{1} 
\lambda^{k}  E_{\beta}(k) & \dots \ \dots & \binom{k}{m-1} 
\lambda^{k}E_{\beta}(k) \\
                         & \ddots& \ddots & \vdots  &\\
                        & & \ddots &  \binom{k}{1} \lambda^{k}  
E_{\beta}(k)  \\[0.2cm]
                        (0) &&  & \lambda^k E_{\beta}(k)   
\label{matrix4}
                \end{pmatrix}.
                \end{equation}
        Now if $A$ is not orthogonal then at least one of the 
following cases will be true.

\medskip

\emph{Case 1:} $A$ has a real eigenvalue $\lambda$ with $|\lambda| 
\neq 1$. Replacing $A$ by $A^{-1}$ if necessary
        we may assume that $| \lambda | > 1$.
        A corresponding block of $A^k$ is an $m \times m$ upper 
diagonal matrix of the form \eqref{matrix4} with $E_{\beta}=1$, and one 
easily checks that
        \[ S = \{ \; s \vecvone : 1 \leq |s| < |\lambda| \;  \} 
\times span(\vecvtwo,\dots,\vecvn)      \]
        is a cross-section.

\medskip

\emph{Case 2:} $A$ has a complex pair of eigenvalues $\lambda e^{\pm 
i\beta}$ with $\lambda \neq 1$, $0 < \beta < \pi$.
        We may again assume that $\lambda>1$.
        A corresponding block of $A^k$  is a $2m \times 2m$
        upper diagonal matrix of form \eqref{matrix4} with 
$E_{\beta}$ a proper rotation.
        Then
                \[      S = \{ \; s \vecvone \lambda^t E_{\beta}(t) : 1 \leq s 
< \lambda^{2\pi/\beta},  \; 0 \leq t < 1 \; \} 
                                \times span(\vecvthree,\dots,\vecvn 
)           \]
        is a cross-section,  which can be checked by using case 2 in Theorem \ref{thm1} and keeping in mind the note
immediately following definition \ref{cross-section}.

\medskip

\emph{Case 3:} $A$ has a real eigenvalue $\lambda = \pm 1$ and at 
least one of the blocks of $A$ belonging to this eigenvalue has 
nontrivial nilpotent part.  Then the corresponding block of $A^k$ is 
of the form
        (\ref{matrix4}) with  $E_{\beta}=1$ and is of at least size 
$2 \times 2$.
       On easily verifies that the set
                \[  S = \left \{ \; s( \vecvone+t \vecvtwo) : s \in 
\R\backslash \{0\}, \; 0 \leq t < 1 \; \right\}  
                        \times span(\vecvthree,\dots,\vecvn  )  \]
        is a cross-section.
\medskip

\emph{Case 4:}  $A$ has a complex pair of eigenvalues 
$e^{\pm  i\beta}$, $ 0 < \beta < \pi$, of modulus one and at least one of 
the blocks of $A$ belonging to this pair has nontrivial nilpotent 
part.        Then the corresponding block of $A^k$ is of
        the form (\ref{matrix4}), with $\lambda =  1$ and $E_{\beta}$ 
a proper rotation,
        and
        \begin{multline*}  S = \Bigl \{ ( p \vecvone + q \vecvthree + 
s \vecvfour) 
                        \begin{pmatrix}
                                E_{\beta}(t)    &  t E_{\beta}(t) \\
                                0               &  E_{\beta}(t)
                        \end{pmatrix}   \\
           : p>0, \ 0 \leq q < \frac{2\pi}{\beta} \, p, 
			\ s \in \R, \ 0 \leq t < 1  \Bigr \} 
                \times  span(\vecvfive,\dots,\vecvn )
        \end{multline*}
         is the desired cross-section, which can be checked by using case 4 in Theorem \ref{thm1} and keeping in mind the note
immediately following definition \ref{cross-section}.  

\medskip

The argument at the end of the proof of theorem \ref{thm1} shows that
if $A$ is orthogonal, then there can not exist a cross-section.  This proves the first assertion.
        
\medskip
The remaining assertions are obvious if $n=1$, or if $n=2$ and $A$ has complex eigenvalues. 
We thus can exclude this situation in what follows,
	so that the cross-section $S$ constructed above has infinite measure.

\medskip

        Next let use prove the second assertion. In order to show 
that $|\det(A)| \neq 1$ is a sufficient condition, we only need to 
distinguish between the first two of the above cases.

        We begin by considering the first case, and we may assume that $|\lambda | >1$. 
Take the cross-section constructed above,
                \[ S = \{ \; s\vecvone  +\vecv : 1 \leq |s| < 
|\lambda|, \  \vecv \in span(\vecvtwo,\dots,\vecvn) \; \},         \]
        partition $ span(\vecvtwo,\dots,\vecvn)$ into a collection 
$\{T_k\}_{k=1}^{\infty}$ of measurable sets of positive, finite 
measure each,
        and set
        \[  S_k =  \{ \; s \vecvone + \vecv : 1 \leq |s| < |\lambda|, 
\  \vecv \in T_k \;  \}, \qquad k=1,2,\dots \]
        Then $\{S_k\}_{k=1}^{\infty}$
        is a partition of $S$ into measurable subsets  of positive, 
finite measure. Pick a collection of positive numbers 
$\{d_k\}_{k=1}^{\infty}$
        so that $\sum_{k=1}^{\infty} d_k = 1$, and pick $n_k \in \Z$ 
such that $\delta^{n_k} \leq \frac{d_k}{m(S_k)}$ where $\delta=|\det(A)|$.
        It follows that
        \[  \tilde S := \bigcup_{k=1}^{\infty} S_kA^{n_k} \] 
        is a cross-section for the discrete action such that
        \[  m(\tilde S) = \sum_{k=1}^{\infty}   \delta^{n_k} 
m(S_k) \leq \sum_{k=1}^{\infty} d_k = 1.  \]

        In the second case, we may assume that $\lambda>1$. 
Start with the above constructed cross-section,
                \[ S = \{ \; s\vecvone \lambda^t  E_{\beta}(t) +\vecv : 1 \leq s 
< \lambda^{2\pi/\beta},\ 0 \leq t < 1, \ 
                        \vecv \in  span(\vecvthree,\dots,\vecvn) \; 
\}  ,\]
        partition $span(\vecvthree,\dots,\vecvn)$ into a collection 
$\{T_k\}_{k=1}^{\infty}$ of 
         measurable subsets of finite, positive measure each, and
        set
        \[  S_k = \{ \; s\vecvone \lambda^t  E_{\beta}(t) + \vecv :  1 \leq s < 
\lambda^{2\pi/\beta},\ 0 \leq t < 1,\ \vecv \in T_k \; \}, 
                                \qquad k=1, 2 \dots \]
        so that $\{S_k\}_{k=1}^{\infty}$ is a partition of $S$ into 
measurable subsets  of positive, finite measure. 
        Continuing as in the first case we have shown sufficiency.
        
        To prove the necessity implication, suppose there exists a 
	cross-section $P$ of finite measure for the
        discrete action and $|\det(A)|=1$. Let $S$ denote the 
	cross-section for the discrete action constructed in part 1 above.
        Then,
        \begin{align*}
         m(P) 	&= \int_{\rnhat}  \chi_P(\vecgamma) \; d \vecgamma 
                                  \ = \ \sum_{i \in \Z} \int_S   \chi_P(\vecgamma A^i) \; d \vecgamma \\
              	&= \sum_{i \in \Z} \int_{\rnhat}  \chi_P(\vecgamma A^i) \, \chi_S(\vecgamma) \; d \vecgamma \\
           	&= \sum_{i \in \Z} \int_{\rnhat}  \chi_P(\vecgamma) \, \chi_S(\vecgamma A^{-i}) \; d \vecgamma \\
                &= \sum_{i \in \Z} \int_P  \chi_S(\vecgamma A^{-i}) \; d \vecgamma \\
                &= \int_{\rnhat} \chi_S(\vecgamma) \; d \vecgamma  \= m(S) \= \infty
        \end{align*}
        which is impossible. Thus, there can not exist a 
cross-section of finite measure. 

\medskip

Finally we will prove the last assertion.    For sufficiency, it is 
enough to assume that all eigenvalues of $A$ have modulus $|\lambda| 
< 1$ so that
                \[  \lim_{ k \to \infty} \|A^k \| = 0. \]
        Choosing each of the above sets $T_k$ to be bounded
	we may assume that the sets $S_k$ are bounded,
	so that there exist integers $n_k$
        such that $S_k A^{n_k}$ is contained in the unit ball. Then 
$\tilde S = \bigcup_{k=1}^{\infty} S_kA^{n_k}$ is the desired bounded
cross-section.

        For necessity, suppose to the contrary that there exists a 
bounded cross-section $\tilde S$, but $A$ has an eigenvalue 
$|\lambda_1|<1$ and an eigenvalue $|\lambda_2|\ge 1$. (The case where 
$|\lambda_1| \leq 1$ and 
        $|\lambda_2| > 1$ is treated similarly). Using the block 
decomposition of $A$
it is easy to see that for almost all $\vecgamma \in \rnhat$, 
either
        \[  \lim_{|k| \to \infty} \|\vecgamma A^k \|  = \infty \]
or, in the special case where no eigenvalue of $A$ lies outside of the unit circle,
        \[  \lim_{k \to - \infty} \|\vecgamma A^k \| = \infty\]
while $\{ \vecgamma A^k: k \ge 0 \}$ is bounded below away from zero.
Thus, for almost all $\vecgamma \in \rnhat  $ there exists a 
constant $M=M(\vecgamma)$ so that
        \[   \|\vecgamma A^k\| > M      \qquad\qquad \forall k \in \Z .
\]
Fix any such $\vecgamma $. Then for sufficiently large scalars $c$, 
the orbit of $c \vecgamma $  does not pass through $\tilde S$,
contradicting the choice of $\tilde S$.
\end{proof}

We note that in the proof of the second assertion,  the sets $T_k$ 
can be chosen
so that the cross-section $\tilde S$ has unit measure.

\begin{remark}  In \cite{LWWW02}, it was obtained as a corollary of their general work that 
\begin{enumerate}
\item  For $A\in GL_n(\R)$ and $D=\{ A^k: k \in \Z \}$, there is a continuous wavelet
if and only if $|\det(A)|\not= 1$.
\item For $A= e^B$ and $D=\{ A^t: t \in \R \}$, there is a continuous wavelet if and only 
if $|\det(A)|\not= 1$.
\end{enumerate}
It is possible to recover these results using the ideas in this section.  
We mention only how to do so in the case that continuous wavelets exist.  Let
$A \in GL_n(\R)$, and let $S$ be a cross-section of Lebesgue measure 1 for the 
discrete action $\vecgamma \to \vecgamma A^k$.  
Then, the function $\psi$ whose Fourier transform equals $\chi_S$ is a continuous wavelet 
for the group $\{ A^k: k \in \Z \}$. If in addition, $A=e^B$, then $\psi$
is also a continuous wavelet for the group $\{ A^t: t \in \R \}$ since
\[	\int_{\R} | \chi_S(\vecgamma A^t)|^2 \; dt \= \int_0^1  \sum_{k \in \Z} | \chi_S(\vecgamma A^tA^k)|^2 \; dt
			\= \int_0^1 1 \; dt \= 1 .	\]
   We note here that the method of proof in \cite{LWWW02}, while
ostensibly constructive, does not easily yield cross-sections of a desirable 
form such as the ones constructed above.  
\end{remark}

\section{Shift-invariant Spaces and Discrete Wavelets}\label{sect13}

Let $A\in GL_n(\R)$ and $\Gamma\subset \R^n$ be a full-rank lattice.  
An $(A, \Gamma)$ orthonormal [resp. Parseval, Bessel]
wavelet \index{wavelet!orthonormal} \index{wavelet!Bessel} of order
$N$ is a collection of functions $\{\psi^i\}_{i=1}^N$ (where here we 
allow the possibility of $N = \infty$) such that  
$$
\{|\det A|^{j/2} \psi^i(A^j \cdot \, + \,  k) : j\in \Z, k\in \Gamma, i=1,\ldots, 
N\}
$$
is an orthonormal basis [resp. Parseval frame, Bessel system] for $L^2(\R^n)$.  
There has been much work done on determining for which
pairs $(A, \Gamma)$ orthonormal wavelets of finite order exist, often with extra desired properties such as
fast decay in time or frequency.

This is not necessary for the proofs that we present.  Of particular importance in determining when orthonormal wavelets exist are the MSF (minimally
supported frequency) wavelets, which are intimately related to wavelet sets.  An $(A, \Gamma)$ multi-wavelet set $K$ of order $L$ is a set that can be
partitioned into subsets $\{K_i\}_{i = 1}^L$ such that $\{\frac{1}{|\det(B)|^{1/2}}\chi_{K_i}\}_{i = 1}^L$ is the Fourier transform of an $(A, \Gamma)$
orthonormal wavelet, where $\Gamma = B \Z^n$, where $B$ is an invertible matrix.  When the order of a multi-wavelet set is 1,
we call it a wavelet set.  These have been studied in detail in
\cite{BMM99, BenLe01, BenLe99, DLS97,  HWW1, HWW2, OS03, W02}. \index{wavelet set}\index{multi-wavelet set}The following fundamental question in this
area remains open, even in the case $L = 1$.

\begin{question} For which pairs $(A, \Gamma)$ and orders $L$ do there exist $(A, \Gamma)$ wavelet sets of order $L$?
\end{question}

It is known that if $A$ is expansive and $\Gamma$ is any full-rank lattice, then there exists an $(A, \Gamma)$ wavelet set of order 1 \cite{DLS97}.  One
can also modify the construction to obtain $(A, \Gamma)$ wavelet sets of any finite order along the lines in Theorem \ref{maintheorem} below.  Diagonal
matrices $A$ for which there exist $(A, \Z^n)$ multi-wavelet sets of finite order were characterized in \cite{S03}.  Theorem \ref{thm2}, part \ref{part2}
above implies that, in order for an $(A, \Gamma)$ multi-wavelet set of finite order to exist, it is necessary that $A$ not have determinant one.   There
is currently no good conjecture as to what the condition on $(A, \Gamma)$ should be for wavelet sets to exist.  It is known that $|\det(A)| \not= 1$ is not
sufficient and that all eigenvalues greater than or equal to 1 in modulus is not necessary.

We begin with the following.

\begin{theorem} \label{wavesetchar}\index{wavelet set!characterization} Let $A\in GL_n(\R)$ and $\Gamma\subset \R^n$ be a full-rank lattice with dual
$\Gamma^*$. The  set $K$ is a multi-wavelet set of order $L$ if and only if
\begin{equation} \label{wavesetchar1} \sum_{\gamma\in \Gamma^*} \chi_K(\xi + \gamma) = L \,\,\, a.e. \,\,\, \xi \in \rnhat, \end{equation}
\begin{equation} \label{wavesetchar2} \sum_{j \in \Z} \chi_K(\xi A^j) = 1 \,\,\, a.e. \,\,\, \xi \in \rnhat. \end{equation}
\end{theorem}

\begin{proof}  The forward direction is very similar to the arguments presented in \cite{DLS97}, so we sketch the proof only.  Let $K$ be a multi-wavelet
set of order $L$.  Partition $K$ into $\{K_i\}_{i = 1}^L$ such that $\frac{1}{|\det(B)|^{1/2}}\chi_{K_i}$ is an $(A, \Gamma)$ multi-wavelet of order
$L$.  Then, since $\chi_{K_i} (\xi A^j)$ is orthogonal to $\chi_{K_k} (\xi A^l)$ for each $(i,j) \not= (k,l)$, it follows that $K_i A^j \cap K_k A^l$ is a
null-set when $(i,j) \not= (k,l)$.  Therefore, $\sum_{j \in \Z} \chi_K(\xi A^j) \le 1$ a.e. $\xi \in \rnhat$.  Moreover, since every $L^2$ function can be
written as the combination of functions supported on $\cup_{j = 1}^\infty K A^j$, it follows that $\sum_{j \in \Z} \chi_K(\xi A^j) = 1$ a.e. $\xi \in
\rnhat$, proving
(\ref{wavesetchar2}).  To see (\ref{wavesetchar1}), since $K_i$ is disjoint from $K_j A^k$ for all $(j,k) \not= (i,0)$, it follows that
$\{\frac{1}{|\det(B)|^{1/2}} e^{2\pi i \langle \xi, \gamma \rangle} : \gamma\in \Gamma\}$ must be an orthonormal basis for $L^2(K_i)$.  This implies
(\ref{wavesetchar1}).

For the reverse direction, it is clear that what is needed is to partition $K$ into $\{K_i\}_{i = 1}^L$ so that each $K_i$ satisfies $\sum_{\gamma\in
\Gamma^*} \chi_{K_i}(\xi + \gamma) = 1$ a.e. $\xi \in \rnhat$.  This will follow from repeated application of the following fact.  Given a measurable set
$K$ such that $\sum_{\gamma\in \Gamma^*} \chi_K(\xi + \gamma) \ge 1$ a.e. $\xi \in \rnhat$, there exists a set $U = U(K) \subset K$ such that 
\begin{equation}\label{ueqn} 
\sum_{\gamma\in \Gamma^*} \chi_U(\xi + \gamma) = 1,\,\,\,\,\, a.e.\,\,\,\,\, \xi \in \rnhat.  
\end{equation}

Now, let $\{V_i\}_{i = 1}^\infty$ be a partition of $\rnhat$ consisting of fundamental regions of $\Gamma^*$; that is, the sets $V_i$ satisfy
$\sum_{\gamma\in \Gamma^*} \chi_{V_i} (\xi + \gamma) = 1$ a.e. $\xi \in \rnhat$.  For a set $M \subset \rnhat$ we define $M^t = \cup_{\gamma\in
\Gamma^*} (M + \gamma).$
 Let 
\[ 
L_0 = K.
\]
Let 
\[K_1 = (V_1 \cap L_0) \cup \bigl(U(L_0) \setminus (V_1 \cap L_0)^t \bigr), \]
where $U(L_0)$ is the subset of $L_0$ satisfying (\ref{ueqn}).  Let $L_1 = L_0 \setminus K_1$, and notice that $L_1$ satisfies (\ref{wavesetchar1})
with the right hand side reduced by 1.  
In general, let 
\[K_i = (V_i \cap L_{i - 1})  \cup\bigr(U(L_{i-1}) \setminus (V_i \cap L_{i - 1})^t \bigr), \]
and
\[
L_i = L_{i-1} \setminus K_i.
\]
In the case that $L$ is finite, this procedure will continue for $L$ steps, resulting in a partition of $K$ with the desired properties.  In this case,
the initial partition $\{V_i\}$ was not necessary.  In the case $L = \infty$, since the $V_i$'s partition $\rnhat$, the union of the $K_i$'s will contain
$K$.  Since the $K_i$'s were constructed to be disjoint and to satisfy (\ref{ueqn}), the proof is complete.

\end{proof}

There is also a soft proof of the reverse direction of Theorem \ref{wavesetchar}, that yields slightly less information about wavelets, but provides some
interesting facts about shift-invariant spaces.  Before turning to the applications of Theorem \ref{wavesetchar}, we provide this second proof.

When $L$ is finite, we call an $(A, \Gamma)$ orthonormal wavelet \index{MSF wavelet!combined}\index{combined MSF wavelet} 
$\{\psi^i\}_{i = 1}^L$ an $(A, \Gamma)$ \emph {combined MSF wavelet} if $\cup_{i = 1}^L \supp(\hat
\psi^i)$ has minimal Lebesgue measure.  This terminology was introduced in \cite{B02}, where it was 
shown that the minimal Lebesgue measure is $L$.  It was also shown that if $\{\psi^i\}_{i = 1}^L$ is a combined MSF wavelet, then there is a
multi-wavelet set $K$ of order $L$ such that $K = \cup_{i = 1}^L \supp(\hat \psi^i)$.

When $L = \infty$, it is not clear what the significance is for the union of the supports of $\hat \psi^i$ to have minimal Lebesgue measure.  For this
reason, we adopt the following definition.  An $(A, \Gamma)$ orthonormal wavelet
$\{\psi^i\}_{i = 1}^L$ is an $(A, \Gamma)$ combined MSF wavelet if $K = \cup_{i = 1}^L \supp(\hat\psi^i)$ is a multi-wavelet set of order $L$.  This
definition agrees with the previous definition in the case $L$ is finite.

Let us begin by recalling some of the basic notions of shift-invariant spaces. \index{shift-invariant space}A closed subspace $V\subset L^2(\R^n)$ is
called {\it shift-invariant} if whenever $f\in V$ and $k\in \Z^n$, 
$f(x + k)\in V$.  The shift-invariant space generated
by the collection of functions $\Phi \subset L^2(\R^n)$ is denoted by 
$\calS (\Phi)$ and given by
\[
\overline{\text {span}}\{\phi (x + k): k\in \Z^n, \phi\in \Phi\}.
\]
Given a shift-invariant space, if there exists a finite set $\Phi 
\subset L^2(\R^n)$ such that $V = \calS(\Phi)$, then we
say $V$ is finitely generated.  In the case $\Phi$ can be chosen to be 
a single function, we say $V$ is a principal shift-invariant (PSI) space.  For further basics about shift-invariant 
spaces, we recommend \cite{Bow00,BDR94,H64}.  We will follow closely the development in \cite{Bow00}.

\begin{proposition}\label{abb} The map $\calT : L^2(\R^n) \to L^2(\T^n, 
\ell^2(\Z^n))$ defined by 
\[
\calT f (x) = (\hat f(x + k))_{k\in \Z^n}
\]
is an isometric isomorphism between $L^2(\R^n)$ and  $L^2(\T^n, \ell^2(\Z^n)),$ where $\T^n = \R^n/\Z^n$ is identified with its fundamental domain,
e.g.~$[0, 1)^n$.
\end{proposition}

In what follows, as in Proposition \ref{abb} we will always assume that $\T^n = \R^n/\Z^n$ is identified with $[0,1)^n$.

A {\it range function} is a mapping \index{range function}
\[
J:\T^n \to \{E\subset \ell^2(\Z^n) : E  \text{ is a closed linear subspace}\}.
\]
The function $J$ is measurable if the associated orthogonal 
projections $P(x) : \ell^2(\Z^n) \to J(x)$ are weakly
operator measurable.  With these preliminaries, we can state an 
important theorem in the theory of shift-invariant
spaces, due to Helson \cite{Bow00}.

\begin{theorem}\label{helsonschar} A closed subspace $V \subset L^2(\R^n)$
is shift-invariant if and only if 
	\[
V = \{f\in L^2(\R^n) : \calT f(x) \in J(x) \ for \ a.e. \  x\in \T^n \},
\]
where $J$ is a measurable range function.  The correspondence between 
$V$ and $J$ is one-to-one under the convention that
the range functions are identified if they are equal a.e.  
Furthermore, if $V = \calS(\Phi)$ for some countable
$\Phi\subset L^2(\R^n)$, then 
\[
J(x) = \overline {\rm span} \{ \calT \phi(x): \phi \in \Phi\}.
\]
\end{theorem}

\begin{definition} \index{dimension function}
The {\rm dimension function} of a shift-invariant space $V$ is the 
mapping $\dim{V}: \T^n \to \N \cup \{0,\infty\}$
given by 
\begin{equation}\label{dimdef}
	\dim{V} (x) = {\rm dim}\, J(x),
\end{equation} 
where $J$ is the range function associated with $V$.  
The {\it spectrum} of $V$ is defined by 
$\sigma(V) = \{x\in \T^n: J(x) \not= \{0\}\}.$  
\end{definition}

We are now ready to state the main result from \cite{Bow00} that we 
will need in this paper.

\begin{theorem}\label{marcinstheorem}  Suppose $V$ is a shift-invariant subspace of $L^2(\R^n)$.  Then $V$ can be decomposed as an orthogonal sum
\begin{equation}\label{vdecomp}
	V = \bigoplus_{i\in \N} \calS(\phi_i),
\end{equation}
where $\{\phi_i(x + k): k\in \Z^n\}$ is a Parseval frame for 
$\calS(\phi_i)$ and $\sigma(\calS(\phi_{i + 1})) \subset
\sigma (\calS(\phi_i))$ for all $i\in \N$.  
Moreover, 
$\dim{{\calS(\phi_i)}}(x) = \|\calT \phi_i(x)\| \in \{ 0, 1\}$  for $i \in \N$, and 
\begin{equation}\label{dimvsi}
	\dim{V}(x) = \sum_{i \in \N} \|\calT \phi_i(x) \| \hskip .2in for \ a.e. \ x \in \T^n.
\end{equation}
\end{theorem}

Finally, there is a folk-lore fact about dimension functions that we recall here.  See  Theorem 3.1 in \cite{BS02a} for discussion and references.

\begin{proposition}\label{intuitivedimfun} Suppose $V$ is a shift-invariant space such that there exists a set $\Phi$ such that
\[
\{\phi(\cdot + k): k\in \Z^n, \phi\in \Phi\}
\]
is a Parseval frame for $V$.  Then 
\begin{equation}\label{otherdefofdimfun}
	\dim{V}(\xi) = \sum_{\phi\in \Phi}\sum_{k\in \Z^n} |\hat\phi(\xi + k)|^2.
\end{equation}
\end{proposition}

The following theorem is a relatively easy application of Theorem 
\ref{marcinstheorem}, which was certainly known in
the case $N < \infty$, and probably known to experts in the theory of 
shift-invariant spaces in this full generality.  It
seems to be missing from the literature, so we 
include a proof.    

\begin{theorem}\label{orthexists}  Let $V$ be a shift-invariant 
subspace of $L^2(\R^n)$.  There exists a collection $\Phi = \{ \phi_i \}_{i=1}^N
\subset L^2(\R^n)$ such that 
\[
\{\phi_i(x + k): i\in \{1,\ldots N\}, k\in \Z^n\}
\]
is an orthonormal basis for $V$ if and only if $\dim{V}(x)  = N$ a.e. 
$x\in \T^n$.
\end{theorem}

\begin{proof}  For the forward direction, it suffices to show that if $\{\phi_i(x + k): 
k\in \Z^n, i=1,\ldots,N\}$ is an orthonormal basis
for the (necessarily shift-invariant) space $V$, then $\dim{V}(x) = N$ for 
a.e.~$x\in \T^n$.  It is easy to see that if
$\{f(x + k): k \in \Z^n\}$ is an orthonormal sequence, then 
$\sum_{k\in \Z^n} |\hat f(\xi + k)|^2 = 1$ a.e. 
Thus, by Proposition \ref{intuitivedimfun}, $\dim{V}(x) = N$ a.e.

For the reverse direction, assume $V$ is a shift-invariant space satisfying $\dim{V}(x) = N$ a.e. $x \in \T^n$. 
Let $\{\phi_i\}_{i = 1}^\infty$ be the collection of functions such that 
(\ref{vdecomp}) is satisfied.  Using the facts that
$\sigma(\calS(\phi_{i+1}))\subset
\sigma(\calS(\phi_i))$ for all $i$,
$\sigma(V) = \T^n$ and (\ref{dimvsi}), it follows that 
\begin{equation}\label{casesphi}
	\sigma(\calS(\phi_i)) = \begin{cases} \T^n &i\le N,\cr 
				0&i > N.\cr
				\end{cases}
\end{equation}
By equation \ref{vdecomp}, we have for $1 \le i \le N$, 
\[
\mathrm{dim}_{{\mathcal{S}(\phi_i)}}(\xi) = 1 = \| \mathcal{T} \phi_i(\xi)\|^2 = \sum_{k \in \Z^n} |\hat\phi_i(\xi + k)|^2
\]
Thus, $\{\phi_i(x + k) : k\in
\Z^n\}$ is an orthonormal basis for $\calS(\phi_i)$.  Since the spaces $\calS(\phi_i)$ are orthogonal, 
$\{\phi_i(x + k): i\in \{1,\ldots N\}, k\in \Z^n\}$
is an orthonormal basis for $V$, as desired.

\end{proof}

We include a proof of the following proposition for completeness.

\begin{proposition} \label{dimsuppw}
Let $V = \{f\in L^2(\R^n) : \text{supp}(\hat f) \subset W\}.$  Then, 
$V$ is shift-invariant and $\dim{V}(\xi) = \sum_{k\in
\Z^n} \chi_W(\xi + k) = \#\{k\in \Z^n: \xi + k \in W\}$ a.e.
\end{proposition}

\begin{proof} Clearly, $V$ so defined is shift-invariant.  
Let $\{e_k : k\in \Z^n\}$ be the standard basis for $\ell^2(\Z^n)$, 
and let $\psi_k$ be defined by $\hat\psi_k = \chi_{(\T^n + k) \cap W}$, 
again for $k\in \Z^n$.  It is easy to see that 
$V =\calS(\Psi)$, where $\Psi = \{\psi_k : k\in \Z^n\}$.  Therefore, by 
Theorem \ref{helsonschar}, 
$J(\xi) = \overline {\text{span}} \{\calT \psi_k (\xi): k\in \Z^n\} = \overline 
		{\text{span}} \{e_k: \xi + k \in W\}$.  
The result
then follows from the definition of dimension function in (\ref{dimdef}).
\end{proof}

\begin{corollary} \label{suffexistence1}
Let $A\in GL_n(\R)$, and $K$ be a measurable subset of $\rnhat$.  If
\begin{equation}\label{calderon}
	\sum_{j\in \Z} \chi_K(\xi A^j) = 1 \hskip .1in a.e. \ \xi \text{ in } \rnhat, 
\end{equation}
and
\[
	\sum_{k\in \Z^n} \chi_{K}(\xi + k) = N \hskip .1in a.e.\ \xi \text{ in } \rnhat,
\]
then there is an $(A, \Z^n)$ orthonormal wavelet of 
order $N$ with $\cup_{i = 1}^N \supp(\hat\psi^i) = K$.
\end{corollary}

\begin{proof}
%If $\Psi$ is an $(A, \Z^n)$ combined MSF wavelet of order $N$ supported on $K$, then by definition there is a multi-wavelet set $\{K_i\}_{i = 1}^N$
%such that $K = \cup_{i = 1}^N K_i$.  The result then follows from Theorem \ref{wavesetchar}.

By Proposition \ref{dimsuppw} and Theorem \ref{orthexists}, 
there exists $\Psi= \{\psi^i\}_{i=1}^N$ such that 
$\{ M_k \hat \psi^i:   k \in \Z^n, \ i=1,\dots,N \}$ is
an orthonormal basis for $L^2(K)$, where $M_k$ denotes modulation by $k$.  Thus, by (\ref{calderon}), 
$\Psi$ is an $(A, \Z^n)$ wavelet.   
%Theorem \ref{wavesetchar} implies that $\Psi$ is combined MSF.  
\end{proof}

The main theorem in this section is given in Theorem \ref{maintheorem}.  Before stating this theorem, we give three results that will be useful in its
proof.

\begin{lemma}\label{easyMinkowski} Let $C\subset \R^n$ be a cone with non-empty interior, $\Gamma\subset \R^n$ be a full-rank lattice, and $T \in
\N$.  Then, the cardinality of $C\cap \Gamma \cap \bigl(\R^n \setminus B_T(0)\bigr)$ is infinity.
\end{lemma}

\begin{proof}  Let $l$ be a line through the origin contained in the interior of $C$.  The set $U = \{x\in \R^n : {\rm dist}(x,
l) < \epsilon\}$ is a centrally symmetric convex set, and 
\begin{equation}\label{emeqn}
((C\cap B_T(0)) \setminus U) \; \mathrm{is\; bounded.}
\end{equation}
By Minkowski's theorem \index{Minkowski's Theorem} (see, for example Theorem 1, Chapter 2, Section 7 in \cite{Lek69} and discussion thereafter), the
cardinality of
$U\cap
\Gamma$ is infinity.  Hence, by (\ref{emeqn}), the result follows.  
\end{proof}

The following proposition was proven in the setting of wave packets in $L^2(\R)$ in \cite{CKSpp}.  We sketch
the proof here in our setting of wavelets.  

\begin{proposition} \label{CKStheorem} \index{wavelet!Bessel} Suppose $A\in GL_n(\R)$ has 
the following property: for all $Z\subset \rnhat$ with positive measure and all 
$q \in \N$, there exist $x_1 \ldots, x_q \in \Z$ such that 
\[
m \bigl(\bigcap_{i = 1}^q  Z A^{x_i} \bigr) > 0.
\]
Then, for every non-zero $\psi\in L^2(\R^n)$, $\psi$ is not an $(A, \Z^n)$ Bessel wavelet.  
\end{proposition}

\begin{proof}
Let $\psi \in L^2(\R^n)$, $\psi \not= 0$. 
Then there exists a set $Z\subset \rnhat$ of positive measure such that $|\hat \psi(\xi)| \ge C > 0$ 
for all $\xi \in Z$.  By reducing to a subset, we may assume that there exists a constant
$K > 0$ such that, for every function $f \in L^2(\rnhat)$ with support in $Z$, we have
\begin{equation*}
	\sum_{k\in \Z^n} |\langle f, M_k \hat \psi \rangle|^2 \ge K\|f\|^2.
\end{equation*}  
Since the operator $D f = |\det(A)|^{1/2} f(\cdot A)$ is unitary, for every $j\in \Z$ and
for each function $f \in L^2(\rnhat)$ supported in $A^{-j} (Z)$, we obtain
\begin{equation}\label{eqLemma47}
\sum_{k\in \Z^n} |\langle f, D^j M_k\hat\psi \rangle|^2 \ge K\|f\|^2.
\end{equation}
By hypothesis, there exist $x_1 \ldots, x_q \in \Z$ such that for 
$U:= \bigl(\bigcap_{i = 1}^q  Z A^{x_i} \bigr)$, we have $m (U) > 0.$
This implies
\begin{eqnarray*}
   	\sum_{j\in \Z, k\in \Z^n} |\langle \chi_U, D^j M_k\hat\psi\rangle|^2 
		&\ge & \sum_{i=1}^q \sum_{k\in \Z^n} |\langle \chi_U, D^{x_i} M_k\hat\psi\rangle|^2\\
		&\ge & \sum_{i=1}^q K\|\chi_U\|^2\\
		& = & q K\|\chi_U\|^2.
\end{eqnarray*}
Thus, since $q$ is arbitrary, $\psi$ is not an $(A, \Z^n)$ Bessel wavelet.  
\end{proof}

\begin{theorem}\label{bonfer} (Bonferroni's Inequality) 
If $\{A_i\}_{i=1}^N$ are measurable subsets of the measurable set $B$ and
$k$ is a positive integer such that
$$\sum_{i=1}^N |A_i| > k |B|,$$
then there exist $1 \le i_1 < i_2 < \cdots < i_{k+1} \le N$ such that 
$$\left|\bigcap_{j=1}^k A_{i_j} \right| > 0.$$
\end{theorem}

\begin{theorem}\label{maintheorem} Let $A\in GL_n(\R)$ with real Jordan form $J$.  The following
statements are equivalent.
\begin{enumerate}
\item\label{main0} For every full-rank lattice $\Gamma \subset \R^n$, there exists a $(J, \Gamma)$ orthonormal
wavelet of order $\infty$.
\item\label{main1} There exists an $(A, \Z^n)$ orthonormal wavelet of
order $\infty$.
\item\label{main2} For every full-rank lattice $\Gamma\subset \R^n$, there
exists an $(A, \Gamma)$ orthonormal wavelet of order $\infty$.
\item\label{main3} There exists a (non-zero) $(A, \Z^n)$ Bessel wavelet of order 1.
\item\label{main4} $J$ is not orthogonal.
\item\label{main5} The matrix $A$ is not similar (over $M_n(\mathbb{C})$)
to a unitary matrix.
\end{enumerate}
\end{theorem}

\begin{proof}  Let us begin by summarizing the known results and obvious
implications.  The implication (\ref{main1}) $\implies$ (\ref{main5}) was
proven Theorem
4.2 \cite{ILP98}.  The implications (\ref{main2}) $\implies$
(\ref{main1}) $\implies$ (\ref{main3}) are obvious, and (\ref{main4})
$\iff$ (\ref{main5}) is standard.

\smallskip\noindent (\ref{main4}) $\implies$ (\ref{main0}).  Let $\Gamma$ be a full-rank lattice with convex fundamental
region $Y$ for $\Gamma^*$.  
By Theorem \ref{wavesetchar}, it suffices to show that there  is a measurable cross-section $S$ for the
discrete action $\xi \to \xi J^k$ satisfying (\ref{wavesetchar2}).  
As in Theorems \ref{thm1} and \ref{thm2}, we break the analysis into cases.  

\smallskip\noindent {\it Case 1: }There is an eigenvalue of $J$ not equal to 1 in modulus.  WLOG, we assume there is an
eigenvalue of modulus greater than 1.  In this case,
$J$ can be written as a block diagonal matrix 
\begin{equation}\label{cdiag}
	\begin{pmatrix} J_1&0\cr 0& J_2 \end{pmatrix},
\end{equation}
where $J_1$ is expansive, and we allow the possibility that rank$(J_1) = {\rm rank}(J)$.  
Let $S$ be an open cross-section for the discrete action $\xi \to \xi J_1^k$.  
Partition $S$ into disjoint open subsets $\{S_i: i\in \N\}$.  For each $i$, choose $k_i$
such that there exists $\gamma_i \in \Gamma^*$ such that $S_i A^{k_i} \times \R^{{\rm rank}(J_2)}\supset (Y + \gamma_i)$. 
Then, 
\[
	\cup_{i = 1}^\infty (S_i A^{k_i} \times \R^{{\rm rank}(J_2)})
\]
is a cross-section satisfying (\ref{wavesetchar2}).

\smallskip \noindent {\it Case 2: }All eigenvalues of $J$ have modulus 1.    This means that we are in case 3 or case 4 of
Theorem \ref{thm2}.  We show that in either of these cases, the cross-section exhibited in Theorem \ref{thm2} satisfies
(\ref{wavesetchar2}).   First, note that $S$ in these cases is a cone of infinite measure with a dense, open subset
$S^\circ$.  Let $B\subset S^\circ$ be an open ball bounded away from the origin satisfying $\overline{B} \subset S^\circ$.  
Let $\delta = {\rm diam}(Y)$.  There exists a $T$ such that 
\[
	S_T := \{tb: t \ge T, b \in B\}
\]
satisfies ${\rm dist}(S_T, \R^n \setminus S) > \delta$.  By Lemma \ref{easyMinkowski}, $\Gamma^* \cap S_T$ has infinite
cardinality, and by choice of $\delta$, $Y + \gamma \subset S$ for each $\gamma\in \Gamma^* \cap S_T$.  Therefore, $S$ is a
cross-section satisfying (\ref{wavesetchar2}).

\smallskip\noindent {\it (\ref{main0}) $\implies$ (\ref{main1}) $\implies$ (\ref{main2})}.  This follows from the following
two facts.  First, $\Gamma$ is a full-rank lattice if and only if there is an invertible matrix $B$ such that $\Gamma =
B\Z^n$.  Second, if $B\in GL_n(\R)$, $\Psi$ is an $(A, \Gamma)$ orthonormal wavelet if and only if 
$\Psi_B := \{ \frac{1}{|\det(B)|^{1/2}} \psi(B^{-1}
\cdot) : \psi \in \Psi\}$ is a $(BAB^{-1}, B\Gamma)$ orthonormal wavelet.  Indeed, (\ref{main0}) $\implies$ (\ref{main1}) is then immediate.  

To see (\ref{main1}) $\implies$ (\ref{main2}), recall that (\ref{main1}) $\implies$ (\ref{main5}).  Thus, if (\ref{main1}) is
satisfied, then $J = B^{-1}AB$ is not orthogonal.  Let $\Gamma$ be a full-rank lattice.  There exists a $(J, B^{-1}\Gamma)$
orthonormal wavelet of order $\infty$, so there exists an $(A, \Gamma)$ orthonormal wavelet of order $\infty$.  

\smallskip\noindent {\it (\ref{main3}) $\implies$ (\ref{main5})}.  Suppose that the real Jordan form of $A$ is orthogonal. 
Then, for any bounded set $Z\subset \rnhat$, there exists $M$ such that for every $k\in \Z, z\in Z$, we have 
$\|z A^k\| \le M$.  Furthermore, if $Z$ has positive measure,  then
\[
\sum_{k\in \Z} m(Z A^k \cap B_M(0)) = \infty.
\]
Therefore, by Bonferroni's inequality, for every $q\in \N$, there exist $k_1, \ldots k_q$ such that 
\[
	m\bigl(\cap_{j = 1}^q Z A^{k_j}  \bigl) > 0.
\]
By Proposition \ref{CKStheorem}, this says that for every non-zero $\psi$, $\psi$ is not an $(A, \Z^n)$ Bessel wavelet.
\end{proof}

\end{document}